\documentclass[a4paper]{article}
\usepackage{amsmath}
\usepackage{amssymb}
\usepackage{amsthm}
\usepackage{amsfonts}
\usepackage{authblk}
\usepackage[english]{babel}
\usepackage{float}
\usepackage{textcomp}
\makeatletter
\usepackage{graphicx}
\usepackage{wrapfig}
\usepackage[misc,geometry]{ifsym}
\usepackage[table]{xcolor}
\usepackage[version=3]{mhchem}
\usepackage{tikz}
\usepackage{filecontents}
\usepackage{sidecap} 
\usepackage{natbib}
\usetikzlibrary{positioning}
\def\m{\phantom{-}}
\newcommand{\eg}{\emph{e.g.,\ }}
\newcommand{\ie}{\emph{i.e.,\ }}
\definecolor{kthgrey}{RGB}{227,229,227}
\definecolor{processblue}{RGB}{36,160,216}
\definecolor{kthblue}{RGB}{25,84,166}
\definecolor{redcol}{RGB}{157,16,45}
\newcommand{\norm}[1]{\left\lVert#1\right\rVert}
\newcommand{\abs}[1]{\left\lvert#1\right\rvert}
\providecommand{\keywords}[1]{\textbf{Keywords:} #1}
\begin{document}
\title{Robust generation of elementary flux modes}
\author[1]{Hildur {\AE}sa Oddsd\'{o}ttir \thanks{Corresponding author: haodd@kth.se}}
\author[2]{Erika Hagrot}
\author[2]{V\'{e}ronique Chotteau }
\author[1]{Anders Forsgren}
\affil[1] {Department of Mathematics, Optimization and Systems Theory, 
	KTH Royal Institute of Technology, 
	SE-100 44 Stockholm, Sweden}
\affil[2]{Division of Industrial Biotechnology/Bioprocess Design, KTH Royal Institute of Technology, 
	Albanova Center,
	SE-106 91 Stockholm, Sweden}
\date{Received: \today / }
\maketitle
\begin{abstract}
Elementary flux modes (EFMs) are vectors defined from a metabolic reaction network, giving the connections between substrates and products. 
EFMs-based metabolic flux analysis (MFA) estimates the flux over each EFM from external flux measurements through least-squares data fitting. 
In previous work we presented an optimization method of column generation type that facilitates EFMs-based MFA when the metabolic reaction network is so large that enumerating all EFMs is prohibitive. 
In this work we extend this model by including errors on measurements in a robust optimization framework. 
In the robust optimization problem, the least-squares data fitting is minimized subject to the error on each metabolite being as unfavourable as it can be, within a given interval. In general, inclusion of robustness may make the optimization problem significantly harder. However, we show that in our case the robust problem can be stated as a convex quadratic programming problem, i.e., of the same form as the original non-robust problem. 
Additionally, we demonstrate that the column-generation technique of the non-robust problem can be extended also to the robust problem. 
Furthermore, the option to indicate intervals on metabolites that are not measured is introduced in this column generation framework. 
The effect of including robustness in the model is evaluated in a case-study, which indicated that the solutions of our non-robust problems are in fact near-optimal also when robustness is considered. On the other hand, the addition of intervals on unmeasured metabolites resulted in a change of optimal solution. Implying that the inclusion of intervals on unmeasured metabolites is more important than the direct consideration of measurement errors, this despite up to $20\%$ errors.
\end{abstract}
\section{Introduction}\let\thefootnote\relax\footnote{\emph{Abbreviations:} MFA, metabolic flux analysis; EFMs, elementary flux modes; CHO, Chinese hamster ovary; Lac, lactate; Glc, glucose}
In previous work we presented a column generation based algorithm for solving the EFMs-based metabolic flux analysis (MFA) problem \citep{Oddsdottir2014}. 
In this work we present a more refined model where the column generation algorithm is combined with robustness. For the sake of completeness a short description of the background follows. A more detailed background can be found in \eg \citet{Oddsdottir2014}.   
A metabolic reaction network is represented by the stoichiometric matrix $A$, which together with the flux vector ($v$) gives the overall change in concentration of each metabolite ($C$).
The rows of the stoichiometric matrix ($A$) refer to either external metabolites ($A_x$) or internal ($A_i$). The flux space is given by a set of vectors $v$ that satisfy the pseudo-steady state assumption and flow direction assumption, 
\begin{equation}\label{eq:Cone}
\left\{v: \begin{bmatrix}
\m A_{i}\\
-A_{i}\\
-I_j
\end{bmatrix}v \leq \begin{bmatrix} 0\\0\\0
\end{bmatrix}, \; j\in J_{\text{irrev}} \right\},
\end{equation}
where $I_j$ is a reduced identity matrix with ones only when $j\in J_{\text{irrev}}$ and   $J_{\text{irrev}}$ is the set of irreversible reactions.
When  all reactions in the network are irreversible \eqref{eq:Cone} is a cone where any ray can be written as a  non-negative linear combination of the extreme rays \citep[Part I.4 Theorem 4.8]{Nemhauser1999}.

 EFMs contain information how extracellular metabolites are connected by detailing which reactions are required for their uptake or production \citep{Llaneras2010}. They are vectors in the flux space, each EFM includes only a minimial set of reactions and is nondecomposable \citep{Klamt2002a}. Further, any vector in the flux space can be denoted as a non-negative linear combination of the EFMs  \citep{Schilling1999,Papin2003},
\begin{equation}
v= \sum_{l=1}^L w_l e_l =Ew, \qquad \gamma \geq 0,
\label{eq:combin}
\end{equation} 
where $e$ denotes a single EFM and the matrix $E$ contains the EFMs as columns. In this sense the EFMs generate the flux space and 
are related to the definition of extreme rays in the cone \eqref{eq:Cone} with only irreversible reactions. In fact when a metabolic network only has irreversible reactions the EFMs and the extreme rays of the cone \eqref{eq:Cone} are equal \citep{Gagneur2004}. We assume, without loss of generality, that the metabolic network has only irreversible reactions, \ie $v_j\geq 0 \; \forall j$.
When the network includes reversible reactions finding all the EFMs is equivalent to finding all the extreme rays of a cone in an extended space where all reactions are irreversible \citep{Gagneur2004,Urbanczik2005}. 
For modest-sized networks enumeration of EFMs is possible and computer programs exist for that purpose, \eg Metatool \citep{VonKamp2006}. However, with increased network size enumeration of EFMs becomes prohibitive \citep{Klamt2002a}. Thus focus has shifted to identify only a subset of the EFMs \citep{DeFigueiredo2009,Kaleta2009,Tabe-Bordbar2013}. 

This work considers the solution of the EFMs-based metabolic flux analysis (MFA) problem \citep[Chapter 5.2]{Provost2006} when the network is large and there are known bounds on measurement errors. EFMs-based MFA uses the decomposition of $v$ given by \eqref{eq:combin} to create a macroscopic network ($A_xE$). The macroscopic fluxes ($w$) are then adjusted so that the flux in the network fit the cell specific external flux measurements ($Q$), \ie 
\begin{equation}\label{eq:EFMsbMFA}
\begin{aligned}
\underset{w}{\text{minimize}} \quad & \frac{1}{2} \|Q-\mathcal{I}A_xE w \|_2^2\\
\text{subject to } \quad 
& w \geq 0.
\end{aligned}
\end{equation}
The formulation given by \eqref{eq:EFMsbMFA} includes multiple repetitions of the same experiments, \ie if $q_{k}$ are results from one repetition, $k$, then $Q^T= [q_{1}^T, \ldots q_{d}^T]$, where $d$ denotes the number of repetitions.
$\mathcal{I}$ is a stacked identity matrix consisting of $d$ identity matrices of size $M_{ext}$ (number of external metabolites) or $\mathcal{I}=[I_{M_{ext}}, \ldots, I_{M_{ext}}]^T$, where $I_{M_{ext}}$ is repeated $d$ times.

EFMs-based MFA as given by \eqref{eq:EFMsbMFA} requires the whole set of EFMs, limiting the application to simplified networks. Methods that can solve the EFMs-based MFA problem without enumerating EFMs exist. One method identifies EFMs beforehand through a series of linear programming (LP) problems \citep{Jungers2011}.
This method is based on the existence of a feasible flux vector $v$, an assumption we will examine in Section \ref{sec:Feasab}. 
In our previous work we introduced a more integrated approach that enables identification of EFMs in conjunction with solving the EFMs-based MFA problem \citep{Oddsdottir2014}. The approach was based on an optimization technique named column generation \citep{Lubbecke2005}, in which large networks can be handled by relying on a master problem and a subproblem that are solved iteratively.
The subproblem gives the master problem a new column every iteration until the solution of the subproblem indicates that the solution of the master problem is optimal to the full optimization problem.

The experimental measurements used to calculate the fluxes in $Q$ in the EFMs-based MFA problem \eqref{eq:EFMsbMFA} are prone to errors, which have been stated to reach at least 20\% \citep{Goudar2009}. For this reason we wanted to consider the sensitivity of the solution with respect to these errors. 
Additionally, in some cases certain metabolites included in the network, are difficult to measure and thus remain unmeasured in the data set. Even though those metabolites are unmeasured in this specific experimental setup some information on their fluxes can be available, and a bound can be added. 
We therefore present an extension to our previous column generation algorithm given by \citet{Oddsdottir2014}. This extension includes both a robust formulation and a version that deals with unmeasured metabolites, while still having the benefit of working with larger networks. In the robust formulation the error on each measurement is assumed bounded, while unmeasured metabolites are given a feasible interval.

In the robust formulation the aim is to minimize the objective function when the assumed errors are such that the objective is as disadvantageous as it can be. For more information on robust optimization please see \citet{Mulvey1995} or \citet{Ben-Tal2009}.
Previous work on robust least-squares mainly focus on errors in both the measurements and the model, in general those formulations are difficult to solve (NP complete) \citep{Ghaoui1997}. However, we show that for this special case, where the errors are only in measurements and bounded by an interval, the robust problem can be formulated as a convex quadratic programming (QP) problem. Furthermore, column generation can be applied to this QP, allowing the problem to be solved without previous enumeration of EFMs.

The paper is outlined as follows. In Section \ref{sec:Feasab} it is shown how the stacked least-squares can be written as least-squares of averages along with an example, showing that metabolic reaction networks do not necessarily have a feasible flux vector for a given set of external measurements. Then we present the main results of this paper in Section \ref{sec:Robustv}; a robust version of the EFMs-based MFA, where column generation can also be applied, along with a version in which intervals for unmeasured metabolites are included. Finally in Section \ref{sec:CS} we present some results comparing the solutions of the robust problem to the EFM-based MFA.

\section{On the Feasibility of the EFMs-based MFA}\label{sec:Feasab}
In this section we examine the uniqueness of the stacked EFMs-based MFA and if there always exists a flux vector that fits the network and measurements exactly. These observations support our main results shown in Section \ref{sec:Robustv}.
To simplify the discussion we consider a problem equivalent to the EFMs-based MFA where a flux vector $v$ is sought, 
\begin{equation}\label{eq:EFMMFA}
\begin{aligned}
\underset{v}{\text{minimize}} \quad & \frac{1}{2} \norm{\mathcal{I}A_xv-Q}_2^2\\
\text{subject to } \quad & A_iv=0,\\
& v \geq 0.
\end{aligned}
\end{equation}
Problem \eqref{eq:EFMMFA} is equivalent to \eqref{eq:EFMsbMFA}, by using the decomposition of $v$ given by \eqref{eq:combin}, thus removing the equality constraint. 
With the stacking of multiple measurements the objective function of \eqref{eq:EFMMFA} seems to represent an overdetermined problem. However, problem \eqref{eq:EFMMFA} can be represented as if it only has one measurement, or as an underdetermined problem by,
\begin{equation*}
 \norm{\mathcal{I}A_{x}v-Q}_2^2=
 \sum_{k=1}^d \norm{ A_{x}v-q_k}_2^2 
=d v^TA_{x}^TA_{x}v -2 \sum_{k=1}^d q_k^TA_{x}v + \sum_{k=1}^d q_k^Tq_k.
\end{equation*}
Thus, the solution $v$ of \eqref{eq:EFMMFA} is equal to the solution of
\begin{equation}\label{eq:EFMMFAav}
\begin{aligned}
\underset{v}{\text{minimize}} \quad & \norm{ A_{x}v-\frac{1}{d}\sum q_k}_2^2\\
\text{subject to } \quad & A_iv=0,\\
& v \geq 0.
\end{aligned}
\end{equation}
Consequently, for a given experimental condition, stacking repetitions is equal to using the average value of the flux measurements. 

In light of that the data fitting can equivalently use the average, \ie only one measurement, it becomes important to consider if there always exists a solution to \eqref{eq:EFMMFAav} with zero residual. 
That is, if 
\begin{equation} \label{eq:exv}
\exists v: \; A_xv=q, \; A_iv=0, \; v \geq 0, \; \text{for any }q.
\end{equation}
For robustness the existence of a solution is especially relevant, because when there is only one measurement that fits the network exactly robust optimization will not give a different solution from the non-robust solution. 
Although, it should be noted that when there are repetitions, or multiple measurements, the solution of the robust optimization can differ from the non-robust solution.

Previous analysis of calculability in networks have considered when there exists a unique $v$ that satisfies \eqref{eq:exv} without the positivity constraint. Hence examining if a network is underdetermined or determined.
In general a full rank matrix has the whole of $\mathcal{R}$ as its range, indicating that there always exists a $v$ such that $A_xv = q$ and $A_iv=0$. When the network is underdetermined this $v$ would not be unique \citep{Klamt2002}. 
However, this assumes that $v$ can be negative in all values. In metabolic networks reactions are often restricted to only one direction. Hence, an underdetermined network may not have a solution for all sets of measurements. A small example of how this can happen follows.
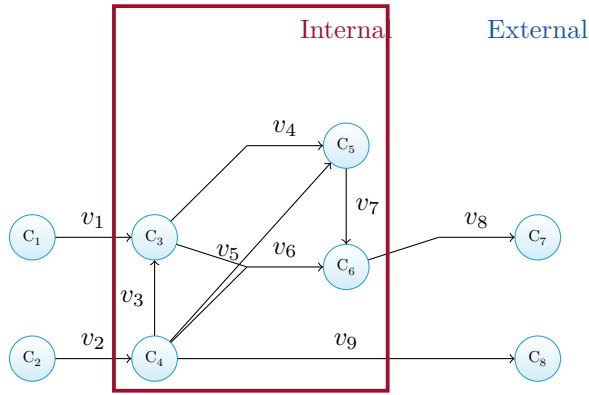
\begin{figure}[H] 
\centering
\begin{tikzpicture}[scale =0.5, state/.style ={ scale=0.6,circle ,top color =white , bottom color = processblue!20 ,draw,processblue , text=black , minimum width =1 cm},interm/.style={fill,circle,inner sep=0pt,scale=0.01}]
\node[interm](cc3){};
\node[state] (c3) [below left =of cc3]{\ce{C3}};
\node[state] (c5) [right =of cc3]{\ce{C5}};
\node[state] (c1) [left =of c3]{\ce{C1}};
\node[state] (c2) [below =of c1]{\ce{C2}};
\node[state] (c4) [right =of c2]{\ce{C4}};
\node[interm] (c5c6) [below right =of c5]{};
\node[state] (c7) [right = of c5c6]{\ce{C7}};
\node[state] (c8) [below = of c7]{\ce{C8}};
\node[interm] (c3c4) [above right =of c4]{};
\node[state] (c6) [right = of c3c4]{\ce{C6}};
\draw[->] (c1) to node[above] {$v_1$} (c3);
\draw[->] (c2) to node[above] {$v_2$} (c4);
\draw[-] (c3) to node {} (cc3);
\draw[->] (cc3) to node[above] {$v_4$} (c5);
\draw[-](c6) to node{} (c5c6);
\draw[-](c3) to node{} (c3c4);
\draw[-](c4) to node{} (c3c4);
\draw[->] (c3c4) to node[above] {$v_6$} (c6);
\draw[->] (c5c6) to node[above] {$v_8$} (c7);
\draw[->] (c4) to node[left]{$v_3$}(c3);
\draw[->] (c4) to node[above]{$v_9$}(c8);
\draw[->] (c5) to node [right] {$v_{7}$}(c6);
\draw[->] (c4) to node [left] {$v_{5}$}(c5);
\draw[line width =0.5mm, color =redcol](-3.5,-6.5) rectangle (3.7,3.7);
\node[redcol] (Int) [above =of c5] {Internal};
\node[kthblue] (Ext) [right =of Int] {External};
\end{tikzpicture}
\caption{A reaction network with underdetermined stoichiometry.}\label{fig:unddRN}
\end{figure}
Consider the network shown in Figure \ref{fig:unddRN}. If one external metabolite is not measured then the network has underdetermined stoichiometry, and thus, there exists a $v$ such that $A_xv=q$ and $A_iv=0$ for any $q$. However depending on which external metabolite is not measured $v$ might not be positive.  
\begin{itemize}
\item If \ce{C1} is not measured, then, depending on what the measurements are, there might not exist a $v \geq 0$ that satisfies the stoichiometry.  With no measurements on \ce{C1}, $v_1$ is free, however $v_2,v_8$ and $v_9$ are fixed from measurements. Flow balance requires that $v_2 \leq v_8+v_9$, additionally if $v_2 \leq v_9$ then the flow through \ce{C4} cannot be fulfilled. Thus, if $v_2$ is too high flow balance can not be fulfilled and no feasible $v$ exists.
\item If \ce{C2} is not measured, then $v_2$ is free, and can be chosen so that the flow to \ce{C7} and \ce{C8} is satisfied, note that any lack of flow from $v_1$ can be compensated by sending through $v_3$.
\end{itemize}
Hence, errors in measurements can lead to nonexistance of a flux vector for the given network that fits the measurements exactly.
\section{The Robust Variant of the EFMs-based MFA}\label{sec:Robustv}
This section contains the main results of this work, here we present an extension of the EFMs-based MFA problem, where errors in $Q$ are taken more directly into consideration. For this purpose we make use of a technique named robust optimization \citep{Mulvey1995,Ben-Tal2009}.

The robust optimization problem is to minimize the residual when the errors in the data give a worst-case scenario outcome, \ie the errors in the data are such that the residual is maximized. Inherent in least-squares is the assumption that the errors are bounded by the two-norm, \ie $\norm{\Delta Q} \leq \beta$. In fact, when the errors are assumed bounded by the two norm, the least squares problem gives the same solution as its robust variant. However, in this work we assume that the errors in $Q$ are bounded by an interval, a more restrictive assumption that might cause the solution to change. The interval is such that $Q_{real}=Q+\Delta Q$ where $\Delta Q_i =[\Delta q_1^T, \ldots \Delta q_d^T]^T$ and $\abs{\Delta q_{ki}} \leq \theta_{ki} \abs{q_{ki}}$, $k$ refers to a specific repetition and $i$ to the metabolite.  
 In order to simplify notation $\theta$ is stacked in the same way as $Q$ and $\Delta Q$, the subindex $s$ then refers to a specific element in those vectors. Note that in general the percentage of error on each metabolite is the same for all repetitions, \ie $\theta_{k_1i}=\theta_{k_2i}$ for all $k_1$ and $k_2$. The robust problem is then given by
\begin{equation}\label{eq:RobOr}
\underset{w \geq 0}{\text{minimize}} \underset{|\Delta Q_{s}| \leq \theta_s |Q_{s}|}{\text{maximize}} \quad  \frac{1}{2} \norm{\mathcal{I} A_{x}E w -Q+\Delta Q}.
\end{equation}
As shown in Appendix \ref{sec:DerRQP}, problem \eqref{eq:RobOr} can equivalently be formulated as a quadratic programming problem in the form 
  \begin{equation}\label{eq:Robdr}
\begin{aligned}
\underset{w,t}{\text{minimize}} \quad & \frac{1}{2} \norm{\mathcal{I}A_{x}E w-Q}^2+\mathbf{1}^T t\\
\text{subject to}\quad & t_s -  \left(\mathcal{I}A_{x}E w-Q\right)_s  \theta_s Q_s\geq0, \quad \forall s\\
& t_s + \left(\mathcal{I}A_{x}E w-Q\right)_s  \theta_s Q_s \geq 0, \quad \forall s\\
&w \geq 0.
\end{aligned}
\end{equation}
When $\theta =0$ the above formulation is equivalent to the EFMs-based MFA \eqref{eq:EFMsbMFA}.

To make the notation more compact we define $\Theta$ and $\tilde{Q}$ as diagonal matrices with $\theta$ and $Q$ on the diagonal, respectively. Further,
the objective function of \eqref{eq:Robdr} can be stated as minimizing the average measure of $q$ over all measurements, 
\begin{equation}\label{eq:Robdrav}
\begin{aligned}
\underset{w,t}{\text{minimize}} \quad & \frac{1}{2} \norm{A_{x}E w-\frac{1}{d}\sum_{k=1}^dq_{k}}^2+\mathbf{1}^T t\\
\text{subject to}\quad & t-\Theta \tilde{Q}\left(\mathcal{I}A_{x}E w-Q\right) \geq 0,\\
& t+\Theta \tilde{Q}\left(\mathcal{I}A_{x}E w-Q\right)  \geq 0,\\
& w\geq0,
\end{aligned}
\end{equation}
The formulation in \eqref{eq:Robdrav} shows that even when the average value gives a zero norm solution of the EFMs-based MFA, the robust solution might be different. 
The reason for this difference can be seen when the constraints in \eqref{eq:Robdrav} are examined. For multiple measurements of the same metabolites $\left(\mathcal{I}A_{x}E w-Q\right)_s$ for each specific measurement will in general not be equal to zero for all $s$, forcing $t$ to increase from zero. 
With enough increase in $t$ the robust solution might deviate from the non-robust solution, \ie increasing $\norm{\mathcal{I}A_{x}E w-Q}$ while decreasing $t$. 
Thereby, giving a non-zero value of $\norm{\mathcal{I}A_{x}E w-Q}$ the robust solution. 
The change in the solution depends on two factors, how far the measurement is from the best least-squares calculated flux and how high the error on that measurement is assumed to be.
No change in the optimal solution is expected when either the measurements are good or the assumed interval is tight, since then $t$ can remain close to zero.
\subsection{Column Generation of the Robust Variant of EFMs-based MFA}
For large networks enumerating all EFMs beforehand is prohibitive. For that reason, we present two problems: A master problem and subproblem that can be solved iteratively to identify the necessary EFMs along with solving problem \eqref{eq:Robdr}. Their derivation can be seen in appendix \ref{sc:MPSP}. The master problem is given by
\begin{subequations} \label{eq:RobMP}
\begin{align}
\underset{w,t}{\text{minimize}} \quad & \frac{1}{2} \norm{\mathcal{I}A_{x}E_B w_B-Q}^2+\mathbf{1}^Tt\\
\text{subject to}\quad & t-\Theta\tilde{Q}\left(\mathcal{I}A_{x}E_Bw_B -Q\right) \geq 0, \label{eq:l1} \\
& t+\Theta \tilde{Q}\left(\mathcal{I}A_{x}E_Bw_B -Q\right) \geq 0, \label{eq:l2}\\
&w_B \geq 0,
\end{align}
\end{subequations}
where the index $B$ indicates that only the known columns of $E$ are used. 
The corresponding subproblem requires information from the solution of the master problem. More specifically the macroscopic fluxes, $w_B$ and the dual solutions, $\lambda_m$ and $\lambda_p$ corresponding to the constraints \eqref{eq:l1} and \eqref{eq:l2} respectively, are required. The subproblem is given by
\begin{equation} \label{eq:RobSubp}
\begin{aligned}
 \underset{e}{\text{minimize}}\quad &\left(\mathcal{I}A_{x}E_Bw_B-Q +\Theta \tilde{Q} \lambda_{m}-\Theta \tilde{Q}\lambda_{p}\right)^T\mathcal{I}A_{x}e \\
 \text{subject to}\quad & A_ie=0,\\
 & \mathbf{1}^Te=1, \\
 & e_j \geq 0 \quad \forall j.
\end{aligned}
\end{equation}
The subproblem \eqref{eq:RobSubp} identifies EFMs \citep{Oddsdottir2014} until the objective function value is non-negative. At that stage the optimal solution of the master problem is also the optimal solution of the full problem. 
\subsection{Inclusion of Intervals on Unmeasured Metabolites in the Robust Variant}
In this section a further extension of the EFMs-based MFA is introduced, where unmeasured metabolites are taken into consideration. Unmeasured metabolites are external metabolites that are a part of the network used but have no measurement data. An example is the metabolite \ce{CO2} a gas that is difficult to measure without special experimental setup. 
Intervals are estimated on those metabolites and modelled with a penalty function. In this way the intervals are allowed to be infeasible for the first few iterations of the column generation. A robust optimization problem that considers feasible intervals on unmeasured metabolites can be stated as,
\begin{equation}\label{eq:Robint}
\begin{aligned}
\underset{w \geq 0}{\text{minimize}}  \;& \left(M_u^T\max(\mathbf{0},(-Q_n^u+A_{x,n}Ew)) \right. \\&+M_l^T\max(\mathbf{0},(Q_n^l-A_{x,n}Ew))  +\\& \;  \left. \underset{\abs{\Delta Q}_k\leq \theta_k \abs{Q_k} }{\text{maximize}} \;  \frac{1}{2} \norm{\mathcal{I}A_{x}E w -Q+\Delta Q} \right),
\end{aligned}
\end{equation}
where $A_{x,n}$ are the rows from the stoichiometric matrix that correspond to the unmeasured metabolites,  $Q_n^u$ and $Q_n^l$ are the upper and lower bounds on the given interval respectively. The quantities $M_u$ and $M_l$ indicate how large the penalty is for not satisfying the specific interval constraint. In general $M_u$ and $M_l$ will be set to a sufficiently large number by the user.
The inner maximization problem is unchanged from \eqref{eq:RobOr} and hence, \eqref{eq:Robint} can be represented as a convex quadratic programming problem,
\begin{subequations}\label{eq:Robintder}
\begin{align}
\underset{w,t,z^u,z^l}{\text{minimize}} \quad & \frac{1}{2} \norm{\mathcal{I}A_{x}E w-Q}^2+\mathbf{1}^T t +M_u^Tz^u + M_l^Tz^l\\
\text{subject to}\quad &  t-\Theta \tilde{Q}\left(\mathcal{I}A_{x}Ew -Q\right) \geq 0,  \\
& t+\Theta \tilde{Q}\left(\mathcal{I}A_{x}Ew -Q\right) \geq 0, \\
&  z^u-A_{x,n}Ew\geq -Q_n^u, \label{eq:lu}\\ 
& z^l+A_{x,n}Ew \geq Q_n^l, \label{eq:ll}\\
& z^u \geq 0,\\
& z^l \geq 0,\\
&w \geq 0.
\end{align}
\end{subequations}
The formulation from \eqref{eq:Robintder} can be solved using column generation where the subproblem generates columns of $E$ by, 
\begin{equation}\label{eq:SPint}
\begin{aligned}
 \underset{e}{\text{minimize}}\quad &\left((\mathcal{I}A_{x}E_Bw_B-Q +\Theta\tilde{Q}\lambda_{m}-\Theta\tilde{Q}\lambda_{p})^T\mathcal{I}A_{x}\right.\\ &\left.+(\lambda_u-  \lambda_l)^TA_{x,n} \right)e\\
 \text{subject to } \quad 
& A_{i}e=0, \\
&\mathbf{1}^Te \leq 1, \\
& e_j \geq 0\; \forall j,
\end{aligned}
\end{equation}
where $\lambda_u$ and $\lambda_l$ are the dual variables corresponding to constraints \eqref{eq:lu} and \eqref{eq:ll}.
\section{Case-Study: Cultivation of CHO Cells} \label{sec:CS}
\subsection{Particulars of the Data}
Data were obtained from the same experimental setup as described in \citet{Oddsdottir2014}. A Chinese hamster ovary (CHO) cell line producing a monoclonal antibody (mAb) was cultivated during 11 days according to a pseudo-perfusion protocol (daily sample collection and medium exchange) to imitate steady-state conditions. The cultivation was carried out in parallel cultures using different medium compositions. Cell-specific metabolic rates (external fluxes) were calculated for the last seven days of culture. Two different media were selected for the present work in order to show extreme situations of our findings, the resulting fluxes are presented in Tables \ref{tb:datam1} (Medium 1) and \ref{tb:datam11} (Medium 5).  In addition to the measured data, an interval on \ce{CO2} flux was estimated as $4.95-7.09$, based on the intervals given by \citet{Goudar2011} and \citet{Aunins1993}.
\begin{table*}[tbh]
\centering
\scriptsize
\rowcolors{1}{white}{kthgrey}
\begin{tabular}{l|rrrrrrr}
\bf{Metabolite} & \bf{$q_{6,1}$} & \bf{$q_{7,1}$} & \bf{$q_{8,1}$} & \bf{$q_{9,1}$} & \bf{$q_{10,1}$} & \bf{$q_{11,1}$} & \bf{$q_{12,1}$}\\ \hline
							Ala & 0.45 & 0.51 & 0.44 & 0.40 & 0.40 & 0.41 & 0.46\\
							Arg &-0.27 &-0.26 &-0.14 &-0.19 &-0.11 &-0.47 &-0.22\\
							Asn &-0.17 &-0.22 &-0.17 &-0.18 &-0.20 &-0.18 &-0.15\\
							Asp &0.04 &0.07 &0.07 &0.06 &0.06 &0.07 &0.07\\
							Biomass &0.61 &0.60 &0.50 &0.70 &0.53 &0.56 &0.65\\
							Cys &-0.09 &-0.11 &-0.07 &-0.09 &-0.09 &-0.07 &-0.05\\
							Glucose (Glc) &-3.52 &-4.06 &-2.64 &-3.26 &-3.96 &-2.92 &-3.43\\
							Gln &-1.60 &-1.97 &-1.61 &-2.38 &-2.31 &-1.90 &-1.71\\
							Glu &0.22 &0.32 &0.25 &0.27 &0.30 &0.33 &0.28\\
							Gly &0.04 &0.07 &0.05 &0.03 &0.03 &0.05 &0.03\\
							His &-0.02 &-0.05 &-0.02 &-0.01 &-0.01 &-0.01 &-0.02\\
							Ile &-0.10 &-0.13 &-0.10 &-0.10 &-0.11 &-0.12 &-0.11\\
							Lactate (Lac) &5.48 &7.40 &5.89 &6.20 &6.78 &7.02 &6.00\\
							Leu &-0.19 &-0.22 &-0.17 &-0.17 &-0.18 &-0.20 &-0.19\\
							Lys &-0.05 &-0.05 &-0.06 &-0.05 &-0.08 &-0.07 &-0.04\\
							Met &-0.05 &-0.07 &-0.04 &-0.05 &-0.06 &-0.04 &-0.03\\
							\ce{NH4+} &1.17 &-- &1.17 &1.15 &1.23 &1.24 &1.18\\
							Phe &-0.10 &-0.12 &-0.10 &-0.09 &-0.09 &-0.12 &-0.12\\
							Pro &-0.10 &-0.14 &-0.09 &-0.11 &-0.11 &-0.13 &-0.10\\
							Ser &-0.00 &0.01 &-0.01 &0.00 &0.00 &-0.03 &0.01\\
							Thr &-0.11 &-0.10 &-0.11 &-0.10 &-0.12 &-0.10 &-0.07\\
							Trp &-0.03 &-0.07 &-0.02 &-0.03 &-0.03 &-0.02 &-0.02\\
							Tyr &-0.09 &-0.12 &-0.08 &-0.08 &-0.08 &-0.12 &-0.10\\
							Val &-0.15 &-0.17 &-0.13 &-0.14 &-0.15 &-0.16 &-0.13\\
							mAb &2.1e-04 &2.3e-04 &1.8e-04 &1.8e-04 &1.7e-04 &2.2e-04 &1.9e-04\\
\end{tabular}
\caption{External fluxes obtained from a CHO cell cultivation, given for each metabolite from the final seven days of the cultivation for medium  1.
The unit is $\text{pmol}\cdot \text{cell}^{-1} \cdot \text{day}^{-1}$, except for Biomass which has the unit $\text{day}^{-1}$.}
\label{tb:datam1}
\end{table*}
\begin{table*}[tbh]
	\centering
	\scriptsize
	\rowcolors{1}{white}{kthgrey}
	\begin{tabular}{l|rrrrrrr}
		\bf{Metabolite} & \bf{$q_{6,5}$} & \bf{$q_{7,5}$} & \bf{$q_{8,5}$} & \bf{$q_{9,5}$} & \bf{$q_{10,5}$} & \bf{$q_{11,5}$} & \bf{$q_{12,5}$}\\ \hline
Ala &0.43 &0.45 &0.46 &0.50 &0.44 &0.40 &0.43\\
Arg &-0.45 &-0.57 &-0.22 &-0.17 &-0.11 &-0.23 &-0.51\\
Asn &-0.18 &-0.22 &-0.20 &-0.22 &-0.21 &-0.19 &-0.17\\
Asp &0.07 &0.08 &0.08 &0.08 &0.09 &0.09 &0.08\\
Biomass &1.11 &0.57 &0.55 &0.59 &0.52 &0.55 &0.56\\
Cys &-0.11 &-0.14 &-0.12 &-0.13 &-0.11 &-0.09 &-0.12\\
Glc &-3.56 &-3.22 &-3.14 &-2.79 &-3.23 &-3.18 &-3.19\\
Gln &-1.79 &-1.71 &-1.60 &-2.41 &-1.81 &-1.85 &-1.77\\
Glu &0.25 &0.31 &0.22 &0.24 &0.22 &0.23 &0.21\\
Gly &-0.00 &-0.06 &0.01 &0.04 &0.03 &0.03 &0.04\\
His &-0.03 &-0.03 &-0.04 &-0.03 &-0.03 &-0.03 &-0.03\\
Ile &-0.13 &-0.20 &-0.11 &-0.11 &-0.11 &-0.11 &-0.13\\
Lac &6.70 &6.02 &5.41 &6.24 &5.61 &5.82 &5.84\\
Leu &-0.21 &-0.30 &-0.19 &-0.19 &-0.18 &-0.19 &-0.22\\
Lys &-0.09 &-0.15 &-0.04 &-0.04 &-0.06 &-0.09 &-0.07\\
Met &-0.05 &-0.07 &-0.06 &-0.05 &-0.06 &-0.04 &-0.04\\
\ce{NH4+} &1.37 &1.45 &1.26 &1.33 &1.31 &1.11 &1.25\\
Phe &-0.12 &-0.11 &-0.10 &-0.11 &-0.09 &-0.07 &-0.13\\
Pro &-0.17 &-0.23 &-0.12 &-0.11 &-0.12 &-0.12 &-0.15\\
Ser &-0.09 &-0.06 &-0.01 &0.03 &0.00 &0.01 &0.01\\
Thr &-0.14 &-0.23 &-0.13 &-0.11 &-0.13 &-0.11 &-0.10\\
Trp &-0.03 &-0.03 &-0.03 &-0.03 &-0.03 &-0.02 &-0.02\\
Tyr &-0.12 &-0.12 &-0.10 &-0.11 &-0.08 &-0.07 &-0.06\\
Val &-0.18 &-0.24 &-0.15 &-0.15 &-0.15 &-0.15 &-0.18\\
mAb &2.3e-04 &2.3e-04 &1.8e-04 &2.2e-04 &1.6e-04 &2.0e-04 &2.1e-04\\
	\end{tabular}
	\caption{External fluxes obtained from a CHO cell cultivation, given for each metabolite from the final seven days of the cultivation for medium 5.
		The unit is $\text{pmol}\cdot \text{cell}^{-1} \cdot \text{day}^{-1}$, except for Biomass which has the unit $\text{day}^{-1}$.}
	\label{tb:datam11}
\end{table*}
\subsubsection{Particulars of the Error on the Data}
The errors on the measurements ($\Delta Q$) are assumed bounded by an error parameter $\theta_s$ that varies for each metabolite but remains constant between repetitions, \ie $\abs{\Delta Q}_s \leq \theta_s \abs{Q}_s$. The estimation of the error parameter was mostly based on the estimated errors of experimental measurements, along with the evaluated variance in the data set. Finally, consistency with the analysis given by \citet{Goudar2009} was ensured.
The values of $\theta$ for each metabolite are given in Table \ref{tb:alpha}.
\begin{table}[tbh]
	\begin{minipage}[c]{0.48\textwidth}
	\centering
	\scriptsize
	\rowcolors{1}{white}{kthgrey}
	\begin{tabular}{l|c}
		\bf{Metabolite} & \bf{Error ($\theta$ [\%])}\\ \hline
Ala & 13.04 \\
Arg & 17.25 \\
Asn & 20.36 \\
Asp & 13.72 \\
Biomass & 17.42 \\
Cys & 17.61 \\
Glc & 14.73 \\
Gln & 15.39 \\
Glu & 13.73 \\
Gly & 15.47 \\
His & 17.10 \\
Ile & 15.31 \\
Lac & 17.52 \\
Leu & 15.49 \\
Lys & 14.55 \\
Met & 13.78 \\
\ce{NH4+} & 13.96 \\
Phe & 16.05 \\
Pro & 15.29 \\
Ser & 15.94 \\
Thr & 15.71 \\
Trp & 15.01 \\
Tyr & 13.58 \\
Val & 23.05 \\
mAb & 18.57 \\
\end{tabular}
\end{minipage}\begin{minipage}[c]{0.5\textwidth}
	\caption{The percentage error on each metabolite, $\theta_i$ (\%)}
	\label{tb:alpha}\end{minipage}
\end{table}
\subsection{Description of the Metabolic Network}
The network used in this study is based on a network available in the literature \citep[Section 2.2]{ZamoranoRiveros2012}. The network was extended in several ways to better fit this study. More reactions were made reversible and some transport reactions were added. The final network consists of 101 reactions, whereof 29 are reversible, and 100 metabolites, whereof 28 are external. Metabolites that are included in the network but not measured are \ce{CO2}, Choline, and Ethanolamine.

In some experiments an external metabolite concentration is set to zero in the medium. This does not exclude the optimal solution from using that metabolite in the optimal solution. Hence, in order to get a solution that fits better with the experimental set-up, columns of $A_{x}$  corresponding to reactions from those metabolites are removed, thus blocking the optimal EFMs from using those reactions. For media 1 and 5 those metabolites are mAb.
\subsection{Technicalities on Normalization}
The results are presented based on a normalized version of the EFMs-based MFA. 
The normalized version aims at fitting the network with the measurements divided by the average value for each specific metabolite in the medium considered. The network is normed similarly by dividing each row of $A_x$ with the average of the measurement for the corresponding metabolite in the medium considered.
 Thus if the average value is defined as,
\begin{equation*}
\bar{q}_{i,g}=\sum_{k=1}^d \frac{q_{i,k,g}}{d},
\end{equation*}
then the external network ($A_x$) and measurements ($Q$) are redefined as follows,
\begin{equation*}
\begin{aligned}
a_{ij}&=\frac{a_{ij}}{\bar{q}_{i,g}} \quad \forall \; j \in J_{ext},\\
q_{i,k,g}&=\frac{q_{i,k,g}}{\bar{q}_{i,g}}.
\end{aligned}
\end{equation*}
Where, $a_{ij}$ is an element from $A_x$ and $J_{ext}$ represents the set of all measured external metabolites in the network. When $\abs{\bar{q}_{i,g}}<0.02$ the value is replaced with $\abs{\bar{q}_{i,g}}=0.02$, in order to avoid dividing by too small values. This minimum is chosen to affect only a few metabolites. For metabolite 1 this affects mAb and Ser, for metabolite 5 this affects those same metabolites along with Gly.
\subsection{Results and Discussion}\label{sec:Res}
\begin{figure}[htb]
	\centering
	\includegraphics[trim=100 0 90 0,clip,width=1\linewidth]{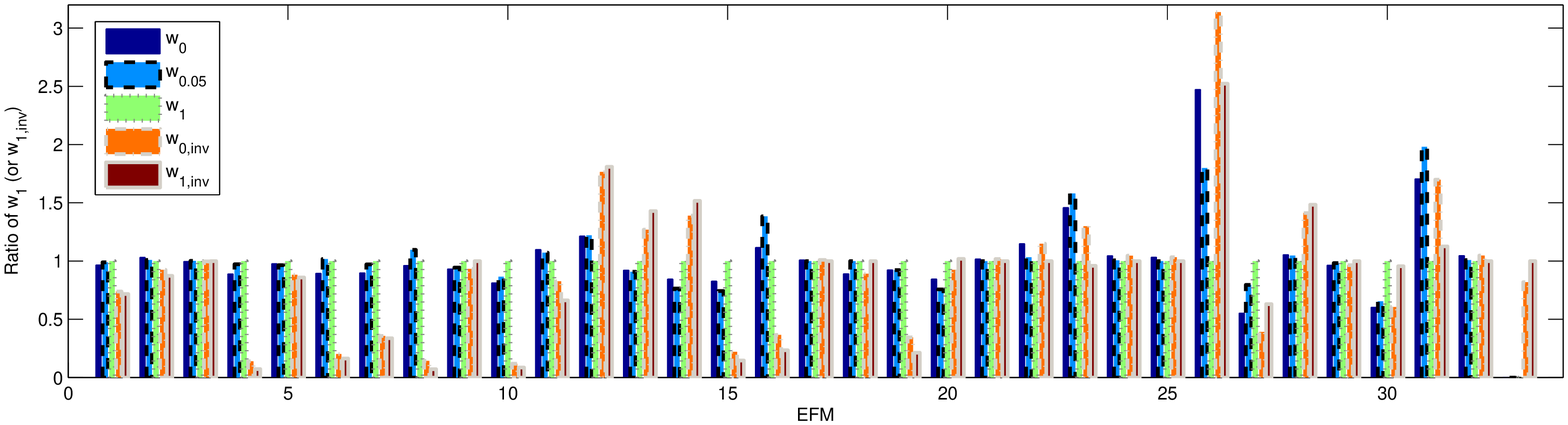}
	\caption{Flux over each EFM, normed with respect to the full error robust solution, for medium 1.}
	\label{fig:M1EFMsw}
\end{figure}

\begin{figure}[htb]
	\centering
	\includegraphics[trim=100 0 90 0,clip,width=1\linewidth]{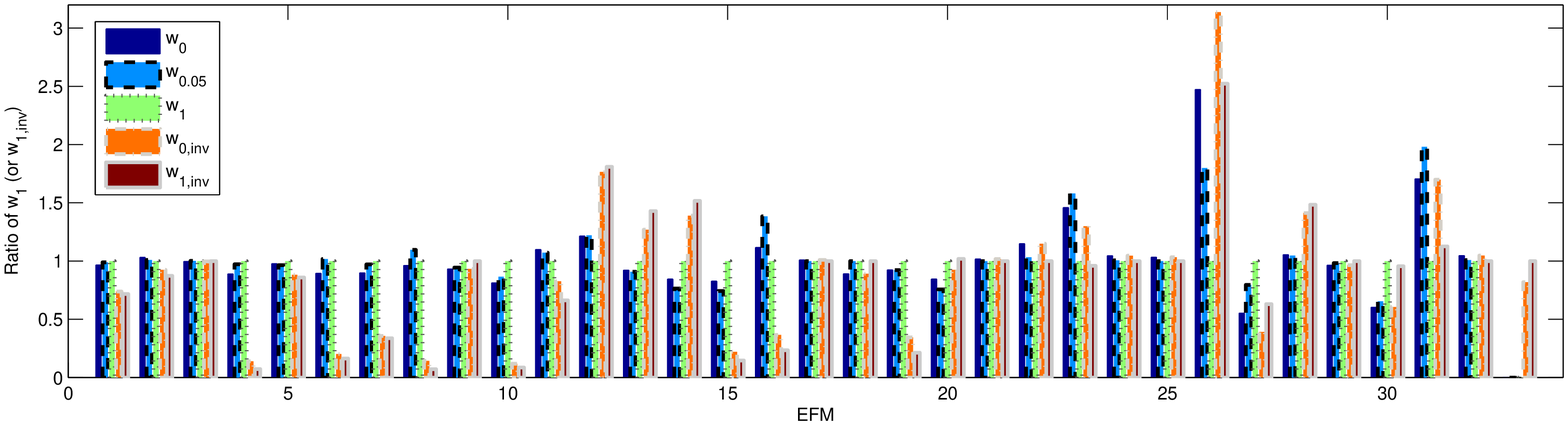}
	\caption{Flux over each EFM, normed with respect to the full error robust solution, for medium 5.}
	\label{fig:M5EFMsw}
\end{figure}

In this section the results for two experimental conditions using two different media are given. In order to demonstrate the difference of the EFMs-based MFA with and without robustness the flux over each EFM and the flux to each external metabolite for three levels of error are shown. The levels of error are 0\%, 5\%, and 100\% of the $\theta$ error given in Table \ref{tb:alpha}. Additionally, the effects of adding an interval are examined by considering the results with a given interval on \ce{CO2} for $\theta$ equal either to zero or 100\%. The 0\% error interval is equivalent to the EFMs-based MFA without robustness. The solution with 100 \% of error interval is referred to as the robust solution. The results are shown in the following tables and figures, where Tables \ref{tb:M1EFMsw} and \ref{tb:M5EFMsw} show the flux over each EFM for medium 1 and 5 respectively. Furthermore, the value of the objective functions for the EFMs-based MFA with and without robustness is shown. 
Figures \ref{fig:M1EFMsw} and \ref{fig:M5EFMsw}, show in the same manner, the flux over each EFM where the flux has been normed with respect to the flux given by the robust solution, this gives an overview of how different the fluxes are for each error interval. 
Tables \ref{tb:M1metw} and \ref{tb:M5metw} 
are similarly constructed but show the flux to each external metabolite. 

\begin{table*}[p]
	\centering
	\scriptsize
	\rowcolors{1}{white}{kthgrey}
	\begin{tabular}{c|p{0.5\linewidth}|rrrrr}
EFM & \bf{Macroscopic Reaction} & $w_0$ & $w_{0.05}$ & $w_1$ & $w_{0,inv}$ & $w_{1,inv}$ \\ \hline
1 &0.5 Glu $\Rightarrow$ 0.5 Ala  + 1 \ce{CO2}  & 5.71 & 5.84 & 5.90 & 4.36 & 4.23\\
2 &0.5 Glc $\Rightarrow$ 1 Lac  & 4.96 & 4.86 & 4.81 & 4.45 & 4.20\\
3 &1 Asn $\Rightarrow$ 1 Lac  & 1.58 & 1.58 & 1.58 & 1.58 & 1.58\\
4 &0.5 Glc  + 0.5 Asn  + 0.5 Ala $\Rightarrow$ 1 Ser  + 0.5 Glu  + 1 \ce{CO2}  & 1.58 & 1.72 & 1.77 & 0.25 & 0.13\\
5 &1 Gln  + 1 Asp $\Rightarrow$ 1 Asn  + 1 Glu  & 1.55 & 1.52 & 1.58 & 1.39 & 1.36\\
6 &1 Ser $\Leftrightarrow$ 1 Gly  & 1.47 & 1.67 & 1.64 & 0.33 & 0.27\\
7 &1 Ala  + 1 \ce{CO2} $\Leftrightarrow$ 1 Asp  & 0.88 & 0.95 & 0.98 & 0.35 & 0.33\\
8 &1 Asp  + 1 Gly $\Rightarrow$ 1 Asn  + 1 \ce{CO2}  & 0.59 & 0.67 & 0.61 & 0.09 & 0.04\\
9 &0.16667 Tyr  + 0.83333 Ala $\Rightarrow$ 1 Asp  & 0.57 & 0.58 & 0.61 & 0.57 & 0.61\\
10 &1 Lac  + 1 Gly $\Rightarrow$ 1 Ala  + 1 \ce{CO2}  & 0.54 & 0.57 & 0.66 & 0.08 & 0.06\\
11 &0.16667 Leu  + 0.16667 Lac  + 0.83333 Ala $\Rightarrow$ 1 Asp  & 0.50 & 0.49 & 0.46 & 0.38 & 0.30\\
12 &1 Gln $\Rightarrow$ 1 Glu  + 1 \ce{NH4+}  & 0.34 & 0.34 & 0.28 & 0.49 & 0.50\\
13 &0.05195 Glc  + 0.0656 Gln  + 0.0278 Ser  + 0.046 Arg  + 0.0552 Thr  + 0.0552 Lys  + 0.0644 Val  + 0.046 Ile  + 0.2596 Leu  + 0.0644 Phe  + 0.0184 Met  + 0.0736 Lac  + 0.0744 Gly  + 0.2032 Glu  + 0.0184 Cys  + 0.0184 His  + 0.046 Pro  + 0.0092 Trp  + 0.006 Ethanolamine  + 0.0171 Choline $\Rightarrow$ 0.084 Asp  + 0.218 \ce{CO2}  + 1 Biomass  & 0.31 & 0.31 & 0.34 & 0.43 & 0.48\\
14 &1 Ala $\Rightarrow$ 1 Lac  + 1 \ce{NH4+}  & 0.29 & 0.26 & 0.34 & 0.47 & 0.52\\
15 &1 Ser $\Rightarrow$ 1 Lac  + 1 \ce{NH4+}  & 0.27 & 0.24 & 0.32 & 0.07 & 0.05\\
16 &1 Gly $\Rightarrow$ 1 \ce{NH4+}  + 1 \ce{CO2}  & 0.23 & 0.29 & 0.21 & 0.08 & 0.05\\
17 &0.5 Val $\Rightarrow$ 0.5 Ala  + 1 \ce{CO2}  & 0.22 & 0.22 & 0.22 & 0.22 & 0.22\\
18 &0.5 Glc  + 1 Arg $\Rightarrow$ 1 Ser  + 1 Glu  & 0.21 & 0.23 & 0.23 & 0.21 & 0.23\\
19 &0.03155 Glc  + 0.0756 Gln  + 0.0368 Asn  + 0.046 Arg  + 0.1982 Thr  + 0.0552 Lys  + 0.0644 Val  + 0.046 Ile  + 0.0828 Leu  + 0.0644 Phe  + 0.0184 Met  + 0.3718 Lac  + 0.1452 Ala  + 0.0164 Glu  + 0.0184 Cys  + 0.0184 His  + 0.046 Pro  + 0.0092 Trp  + 0.006 Ethanolamine  + 0.0171 Choline $\Rightarrow$ 0.3726 \ce{CO2}  + 1 Biomass  & 0.17 & 0.17 & 0.18 & 0.06 & 0.04\\
20 &0.5 Asp  + 0.25 Phe $\Rightarrow$ 0.25 Ala  + 0.5 Glu  + 1 \ce{CO2}  & 0.14 & 0.13 & 0.17 & 0.16 & 0.17\\
21 &0.16667 Lys  + 0.33333 Lac  + 0.66667 Ala $\Rightarrow$ 1 Asp  & 0.14 & 0.14 & 0.14 & 0.14 & 0.14\\
22 &1 Ala  + 0.2 Trp $\Rightarrow$ 0.8 Asp  + 0.4 Glu  & 0.12 & 0.11 & 0.11 & 0.12 & 0.11\\
23 &0.03155 Glc  + 0.0912 Gln  + 0.0686 Ser  + 0.0368 Asn  + 0.046 Arg  + 0.0552 Thr  + 0.0552 Lys  + 0.0644 Val  + 0.046 Ile  + 0.0828 Leu  + 0.3296 Phe  + 0.0184 Met  + 0.0736 Ala  + 0.0744 Gly  + 0.0008 Glu  + 0.0184 Cys  + 0.0184 His  + 0.046 Pro  + 0.0092 Trp  + 0.006 Ethanolamine  + 0.0171 Choline $\Rightarrow$ 0.2092 Asp  + 0.3064 \ce{CO2}  + 1 Biomass  & 0.11 & 0.12 & 0.08 & 0.10 & 0.07\\
24 &1 Pro $\Rightarrow$ 1 Glu  & 8.52e-2 & 8.19e-2 & 8.15e-2 & 8.52e-2 & 8.15e-2\\
25 &1 Ile $\Rightarrow$ 1 Glu  + 1 \ce{CO2}  & 8.25e-2 & 7.99e-2 & 7.99e-2 & 8.25e-2 & 7.99e-2\\
26 &1 Cys $\Rightarrow$ 1 Lac  + 1 \ce{NH4+}  & 5.68e-2 & 4.13e-2 & 2.29e-2 & 7.20e-2 & 5.79e-2\\
27 &1 Gly  + 1 Cys $\Rightarrow$ 1 Asn  & 5.23e-2 & 7.53e-2 & 9.46e-2 & 3.71e-2 & 5.97e-2\\
28 &1 Thr $\Rightarrow$ 0.5 Gly  + 0.5 Glu  + 0.5 \ce{CO2}  & 4.43e-2 & 4.39e-2 & 4.19e-2 & 5.92e-2 & 5.22e-2\\
29 &1 Ser  + 1 Met $\Rightarrow$ 1 Asp  + 1 Cys  & 3.78e-2 & 3.85e-2 & 3.92e-2 & 3.78e-2 & 3.92e-2\\
30 &1 Asp  + 1 His $\Rightarrow$ 1 Asn  + 1 Glu  & 0.56e-2 & 0.61e-2 & 0.93e-2 & 0.56e-2 & 0.89e-2\\
31 &1 His $\Rightarrow$ 1 Glu  + 1 \ce{NH4+}  & 0.56e-2 & 0.65e-2 & 0.33e-2 & 0.56e-2 & 0.37e-2\\
32 & AAs $\Rightarrow$ 1 mAb  & 0.02e-2 & 0.02e-2 & 0.02e-2 & 0.02e-2 & 0.02e-2\\
33 &0.5 Glc $\Rightarrow$ 1 \ce{CO2}  & - & - & - & 1.84 & 2.24\\
Norm & $\norm{A_{x}E w-\frac{1}{d}\sum_{k=1}^dq_{k}}^2$  & 0 & 0.22 & 0.42 & 0 & 0.42 \\
Rob. N.& $\norm{A_{x}E w-\frac{1}{d}\sum_{k=1}^dq_{k}}^2+\mathbf{1}^T t_\theta$ & 463.38 & 444.79 & 442.41 & 463.80 & 442.41
	\end{tabular}
	\caption{Fluxes of each EFM for medium 1, solved with varied assumed error. $w_0$, $w_{0.05}$, and $w_1$ indicate the results with $\theta$ at $0$, $5$, and $100$ of the values given in Table \ref{tb:alpha} \% respectively.  $w_{0,inv}$ and $w_{1,inv}$ indicate the results with an interval given for \ce{CO2} with $\theta$ at $0$ and $100$ \% respectively. $t_\theta=\abs{\tilde{Q}_{\theta}\left(\mathcal{I}A_{x}Ew -Q\right)}$ where $\theta$ is given by Table \ref{tb:alpha}.
		The unit of the fluxes is $\text{pmol}\cdot \text{cell}^{-1} \cdot \text{day}^{-1}$, except for Biomass for which the unit is $\text{day}^{-1}$. AAs denotes amino acids.}
	\label{tb:M1EFMsw}
\end{table*}

\begin{table*}[tbh]
	\begin{minipage}[c]{0.7\textwidth}
	\centering
	\scriptsize
	\rowcolors{1}{white}{kthgrey}
	\begin{tabular}{c|rrrrr}
		\bf{Metabolite} & $q_0$ & $q_{0.05}$ & $q_1$ & $q_{0,inv}$ & $q_{1,inv}$ \\ \hline
Ala & 0.44 & 0.44 & 0.44 & 0.44 & 0.44\\
Arg & -0.24 & -0.26 & -0.26 & -0.24 & -0.26\\
Asn & -0.18 & -0.18 & -0.18 & -0.18 & -0.18\\
Asp & 0.06 & 0.07 & 0.07 & 0.06 & 0.07\\
Biomass & 0.59 & 0.60 & 0.60 & 0.59 & 0.60\\
\ce{CO2} & 8.39 & 8.76 & 8.83 & 6.98 & 7.07\\
Choline & -0.01 & -0.01 & -0.01 & -0.01 & -0.01\\
Cys & -0.08 & -0.09 & -0.09 & -0.08 & -0.09\\
Ethanolamine & -36e-4 & -36e-4 & -36e-4 & -36e-4 & -36e-4\\
Glc & -3.40 & -3.43 & -3.43 & -3.40 & -3.43\\
Gln & -1.93 & -1.91 & -1.90 & -1.93 & -1.90\\
Glu & 0.28 & 0.28 & 0.28 & 0.28 & 0.28\\
Gly & 0.04 & 0.05 & 0.05 & 0.04 & 0.05\\
His & -0.02 & -0.02 & -0.02 & -0.02 & -0.02\\
Ile & -0.11 & -0.11 & -0.11 & -0.11 & -0.11\\
Lac & 6.40 & 6.20 & 6.20 & 6.40 & 6.20\\
Leu & -0.19 & -0.19 & -0.19 & -0.19 & -0.19\\
Lys & -0.06 & -0.06 & -0.06 & -0.06 & -0.06\\
Met & -0.05 & -0.05 & -0.05 & -0.05 & -0.05\\
NH4 & 1.19 & 1.18 & 1.18 & 1.19 & 1.18\\
Phe & -0.10 & -0.10 & -0.10 & -0.10 & -0.10\\
Pro & -0.11 & -0.11 & -0.11 & -0.11 & -0.11\\
Ser & -0.00 & -0.00 & -0.01 & -0.00 & -0.01\\
Thr & -0.10 & -0.10 & -0.10 & -0.10 & -0.10\\
Trp & -0.03 & -0.03 & -0.03 & -0.03 & -0.03\\
Tyr & -0.10 & -0.10 & -0.10 & -0.10 & -0.10\\
Val & -0.15 & -0.15 & -0.15 & -0.15 & -0.15\\
mAb & 2e-4 & 2e-4 & 2e-4 & 2e-4 & 2e-4
	\end{tabular}
\end{minipage}\begin{minipage}[c]{0.28\textwidth}
	\caption{Fluxes of each metabolite for medium 1, solved with varied assumed error. $w_0$, $w_{0.05}$, and $w_1$ indicate the results with $\theta$ at $0$, $5$, and $100$ \% of its given value in Table \ref{tb:alpha} respectively.  $w_{0,inv}$ and $w_{1,inv}$ indicate the results with an interval given for \ce{CO2} with $\theta$ at $0$ and $100$ \% respectively.
		The unit of the fluxes is $\text{pmol}\cdot \text{cell}^{-1} \cdot \text{day}^{-1}$, except for Biomass for which the unit is $\text{day}^{-1}$.}
	\label{tb:M1metw}\end{minipage}
\end{table*}
Considering the results for medium one shown in Table \ref{tb:M1EFMsw} it can be noted that introducing robustness has no impact on which EFMs are required. Further, the change in the norm and robust norm is low, indicating that the impact of robustness on the solution is low. When the flux to each metabolite shown in Table \ref{tb:M1metw} is considered, there are some minor adjustments in the solution when robustness is added, however the larges change comes when an interval on \ce{CO2} is enforced, resulting in a lower flow to \ce{CO2}.

Considering the results from medium five shown in Table \ref{tb:M5EFMsw} it can be noted that for all levels of robustness one more EFM is required in the solution in order to reach optimality. When the flux to \ce{CO2} is bounded a few EFMs are dropped and new ones are introduced. 
In general introducing a new metabolite in the model by giving an interval of its value, \ie without experimental measurements, has a large impact on the solution than introducing robustness. Furthermore the intervals can be satisfied without a large change in the norm of the data fitting. This indicates that giving realistic intervals on unmeasured metabolites can help guide the solution. These results show that, here, larger improvements of the residual error are brought by improving the model structure, i.e. introducing new metabolites, compared to taking into account the robust solution, which address the parameter estimation.

\begin{table*}[p]
	\centering
	\tiny
	\rowcolors{1}{white}{kthgrey}
	\begin{tabular}{c|p{0.6\linewidth}|rrrrr}
EFM & \bf{Macroscopic Reaction} & $w_0$ & $w_{0.05}$ & $w_1$ & $w_{0,inv}$ & $w_{1,inv}$ \\ \hline
1 &0.5 Glc $\Rightarrow$ 1 Lac  & 3.20 & 3.82 & 2.30 & 2.69 & 2.02\\
2 &0.5 Glu $\Rightarrow$ 0.5 Ala  + 1 \ce{CO2}  & 2.64 & 2.72 & 3.13 & 2.53 & 2.48\\
3 &0.5 Glc  + 1 Ala $\Rightarrow$ 1 Ser  + 1 Lac  & 2.62 & 1.93 & 3.44 & 3.22 & 3.72\\
4 &0.5 Gln  + 0.5 Lac $\Rightarrow$ 1 Asp  & 2.05 & 2.27 & 2.20 & 1.49 & 1.62\\
5 &1 Ser $\Leftrightarrow$ 1 Gly  & 1.94 & 2.15 & 2.33 & 2.45 & 2.32\\
6 &1 Asp  + 1 Gly $\Rightarrow$ 1 Asn  + 1 \ce{CO2}  & 1.89 & 2.09 & 2.26 & 1.12 & 1.29\\
7 &1 Asn $\Rightarrow$ 1 Lac  & 1.60 & 1.70 & 1.74 & 1.60 & 1.69\\
8 &0.5 Glu $\Rightarrow$ 0.5 Lac  + 0.5 \ce{NH4+}  + 1 \ce{CO2}  & 1.43 & 1.62 & 1.09 & 1.13 & 1.27\\
9 &0.5 Asn  + 1 Lac $\Rightarrow$ 0.5 Ala  + 0.5 Glu  + 1 \ce{CO2}  & 1.24 & 1.32 & 1.11 & 3e-7 & --\\
10 &1 Ser $\Rightarrow$ 1 Ala  & 1.01 & 0.22 & 1.60 & 1.07 & 1.89\\
11 &1 Gln  + 1 Asp $\Rightarrow$ 1 Asn  + 1 Glu  & 0.54 & 0.35 & 0.37 & 0.73 & 0.69\\
12 &0.5 Glc  + 1 Arg $\Rightarrow$ 1 Ser  + 1 Glu  & 0.36 & 0.46 & 0.49 & 0.35 & 0.50\\
13 &0.25 Tyr $\Rightarrow$ 0.25 Glu  + 1 \ce{CO2}  & 0.34 & 0.26 & 0.28 & 0.26 & 0.29\\
14 &1 Asn $\Rightarrow$ 1 Asp  + 1 \ce{NH4+}  & 0.32 & 0.20 & 0.46 & 0.38 & 0.42\\
15 &0.25465 Glc  + 0.0656 Gln  + 0.0468 Asn  + 0.0616 Asp  + 0.046 Arg  + 0.0552 Thr  + 0.0552 Lys  + 0.0644 Val  + 0.046 Ile  + 0.5446 Leu  + 0.0644 Phe  + 0.0184 Met  + 0.0736 Ala  + 0.0744 Gly  + 0.0184 Cys  + 0.0184 His  + 0.046 Pro  + 0.0092 Trp  + 0.006 Ethanolamine  + 0.0171 Choline $\Rightarrow$ 0.3668 Glu  + 0.9804 \ce{CO2}  + 1 Biomass  & 0.30 & 0.27 & 0.23 & 2e-5 & --\\
16 &0.027502 Glc  + 0.063938 Gln  + 0.066862 Ser  + 0.035867 Asn  + 0.044834 Arg  + 0.053801 Thr  + 0.053801 Lys  + 0.32125 Val  + 0.044834 Ile  + 0.080702 Leu  + 0.32125 Phe  + 0.017934 Met  + 0.072515 Gly  + 0.025731 Glu  + 0.017934 Cys  + 0.017934 His  + 0.044834 Pro  + 0.0089669 Trp  + 0.005848 Ethanolamine  + 0.016667 Choline $\Rightarrow$ 0.024951 Lac  + 0.36569 Ala  + 1 \ce{CO2}  + 0.97466 Biomass  & 0.26 & 0.08 & 0.09 & 0.04 & 0.03\\
17 &1 Gln $\Rightarrow$ 1 Glu  + 1 \ce{NH4+}  & 0.25 & 0.27 & 0.28 & 0.32 & 0.23\\
18 &0.52941 Val  + 0.17647 Leu  + 0.058824 Pro $\Rightarrow$ 0.058824 Arg  + 0.52941 Glu  + 1 \ce{CO2}  & 0.12 & 0.22 & 0.22 & 0.23 & 0.20\\
19 &1 Thr $\Rightarrow$ 1 Asp  & 0.10 & 0.10 & 0.10 & 0.10 & 0.10\\
20 &0.33333 Asn  + 0.33333 Asp  + 0.33333 Lys  + 0.33333 Pro $\Rightarrow$ 0.33333 Arg  + 0.66667 Glu  + 1 \ce{CO2}  & 0.10 & 0.12 & 0.12 & 0.09 & 0.12\\
21 &1 Ile  + 1 Cys $\Rightarrow$ 1 Asp  + 1 Glu  & 0.09 & 0.08 & 0.07 & 0.09 & 0.08\\
22 &0.16667 Glc  + 0.33333 Asn  + 0.33333 Met  + 0.33333 Pro $\Rightarrow$ 0.33333 Arg  + 1 Lac  + 0.33333 \ce{CO2}  & 0.08 & 0.05 & 0.05 & 0.06 & 0.05\\
23 &0.25 Ala  + 0.25 Trp $\Rightarrow$ 0.5 Glu  + 1 \ce{CO2}  & 0.07 & 0.07 & 0.06 & 0.06 & 0.08\\
24 &1 Pro $\Rightarrow$ 1 Glu  & 0.05 & 0.05 & 0.05 & 0.05 & 0.05\\
25 &0.028217 Glc  + 0.0756 Gln  + 0.0686 Ser  + 0.0368 Asn  + 0.0616 Arg  + 0.2044 Tyr  + 0.0552 Thr  + 0.0552 Lys  + 0.0644 Val  + 0.2228 Ile  + 0.0828 Leu  + 0.0368 Phe  + 0.0184 Met  + 0.0744 Gly  + 0.0008 Glu  + 0.0184 Cys  + 0.0184 His  + 0.046 Pro  + 0.0092 Trp  + 0.006 Ethanolamine  + 0.0171 Choline $\Rightarrow$ 0.224 Ala  + 0.4956 \ce{CO2}  + 1 Biomass  & 0.04 & 0.15 & 0.17 & 0.07 & 0.05\\
26 &0.028217 Glc  + 0.0656 Gln  + 0.5834 Ser  + 0.046 Arg  + 0.0276 Tyr  + 0.0552 Thr  + 0.3204 Lys  + 0.0644 Val  + 0.046 Ile  + 0.0828 Leu  + 0.0368 Phe  + 0.5332 Met  + 0.0744 Gly  + 0.0264 Glu  + 0.0184 His  + 0.046 Pro  + 0.0092 Trp  + 0.006 Ethanolamine  + 0.0171 Choline $\Rightarrow$ 0.3484 Ala  + 0.468 \ce{NH4+}  + 0.9736 \ce{CO2}  + 0.4964 Cys  + 1 Biomass  & 0.03 & 0.05 & 0.06 & 0.04 & 0.05\\
27 &1 Cys $\Rightarrow$ 1 Ala  & 0.03 & 0.06 & 0.07 & 0.04 & 0.06\\
28 &0.024997 Glc  + 0.058115 Gln  + 0.060773 Ser  + 0.032601 Asn  + 0.06343 Asp  + 0.040751 Arg  + 0.024451 Tyr  + 0.048901 Thr  + 0.048901 Lys  + 0.057052 Val  + 0.040751 Ile  + 0.073352 Leu  + 0.032601 Phe  + 0.0163 Met  + 0.065911 Gly  + 0.038979 Cys  + 0.0163 His  + 0.050319 Pro  + 0.25 Trp  + 0.0053154 Ethanolamine  + 0.015149 Choline $\Rightarrow$ 0.022679 Lac  + 0.17665 Ala  + 1 \ce{CO2}  + 0.8859 Biomass  & 0.02 & 0.03 & 0.04 & 0.03 & 0.02\\
29 &1 Lac  + 1 His $\Rightarrow$ 1 Ala  + 1 Glu  & 99e-4 & 118e-4 & 175e-4 & 100e-4 & 151e-4\\
30 &1 His $\Rightarrow$ 1 Glu  + 1 \ce{NH4+}  & 99e-4 & 92e-4 & 35e-4 & 99e-4 & 59e-4\\
31 & AAs $\Rightarrow$ 1 mAb  & 2e-4 & 2e-4 & 2e-4 & 2e-4 & 2e-4\\
32 &0.33333 Asp  + 0.33333 Phe $\Rightarrow$ 0.66667 Glu  + 1 \ce{CO2}  & - & 0.18 & 0.18 & 0.21 & 0.24\\
33 &1 Lac  + 1 Gly $\Rightarrow$ 1 Ala  + 1 \ce{CO2}  & - & - & - & 1.28 & 0.95\\
34 &0.028217 Glc  + 0.0656 Gln  + 0.0686 Ser  + 0.0368 Asn  + 0.0716 Asp  + 0.046 Arg  + 0.0276 Tyr  + 0.0552 Thr  + 0.0552 Lys  + 0.0644 Val  + 0.046 Ile  + 0.343 Leu  + 0.0368 Phe  + 0.0184 Met  + 0.2346 Lac  + 0.0736 Ala  + 0.0744 Gly  + 0.0184 Cys  + 0.0184 His  + 0.046 Pro  + 0.0092 Trp  + 0.006 Ethanolamine  + 0.0171 Choline $\Rightarrow$ 0.2238 Glu  + 0.0368 \ce{CO2}  + 1 Biomass  & - & - & - & 0.45 & 0.42\\
& $\norm{A_{x}E w-\frac{1}{d}\sum_{k=1}^dq_{k}}^2$  & 0.00 & 0.69 & 2.27 & 0.01 & 2.27 \\
& $\norm{A_{x}E w-\frac{1}{d}\sum_{k=1}^dq_{k}}^2+\mathbf{1}^T t_\theta$ &1151.23 & 1073.57 & 941.05 & 1154.05 & 945.09
	\end{tabular}
	\caption{Fluxes of each EFM for medium 5, solved with varied assumed error. $w_0$, $w_{0.05}$, and $w_1$ indicate the results with $\theta$ at $0$, $5$, and $100$ of the values given in Table \ref{tb:alpha} \% respectively.  $w_{0,inv}$ and $w_{1,inv}$ indicate the results with an interval given for \ce{CO2} with $\theta$ at $0$ and $100$ \% respectively. $t_\theta=\abs{\tilde{Q}_{\theta}\left(\mathcal{I}A_{x}Ew -Q\right)}$ where $\theta$ is given by Table \ref{tb:alpha}.
		The unit of the fluxes is $\text{pmol}\cdot \text{cell}^{-1} \cdot \text{day}^{-1}$, except for Biomass for which the unit is $\text{day}^{-1}$. AAs denotes amino acids.}
	\label{tb:M5EFMsw}
\end{table*}

\begin{table*}[tbh]
	\begin{minipage}[c]{0.7\textwidth}
	\centering
	\scriptsize
	\rowcolors{1}{white}{kthgrey}
	\begin{tabular}{c|rrrrr}
		\bf{Metabolite} & $q_0$ & $q_{0.05}$ & $q_1$ & $q_{0,inv}$ & $q_{1,inv}$\\ \hline
Ala & 0.44 & 0.44 & 0.44 & 0.45 & 0.44\\
Arg & -0.32 & -0.42 & -0.45 & -0.32 & -0.45\\
Asn & -0.20 & -0.20 & -0.20 & -0.20 & -0.20\\
Asp & 0.08 & 0.08 & 0.08 & 0.08 & 0.08\\
Biomass & 0.64 & 0.57 & 0.57 & 0.64 & 0.57\\
\ce{CO2} & 8.47 & 9.11 & 8.95 & 7.09 & 7.09\\
Choline & -0.01 & -0.01 & -0.01 & -0.01 & -0.01\\
Cys & -0.12 & -0.12 & -0.12 & -0.12 & -0.12\\
Ethanolamine & -38e-4 & -34e-4 & -34e-4 & -38e-4 & -34e-4\\
Glc & -3.19 & -3.19 & -3.19 & -3.16 & -3.14\\
Gln & -1.85 & -1.79 & -1.79 & -1.84 & -1.77\\
Glu & 0.24 & 0.23 & 0.23 & 0.24 & 0.24\\
Gly & 0.01 & 0.02 & 0.03 & 0.01 & 0.03\\
His & -0.03 & -0.03 & -0.03 & -0.03 & -0.03\\
Ile & -0.13 & -0.13 & -0.13 & -0.13 & -0.11\\
Lac & 5.95 & 5.84 & 5.84 & 6.00 & 6.24\\
Leu & -0.21 & -0.21 & -0.19 & -0.21 & -0.19\\
Lys & -0.07 & -0.09 & -0.09 & -0.07 & -0.09\\
Met & -0.05 & -0.05 & -0.05 & -0.05 & -0.05\\
\ce{NH4+} & 1.30 & 1.31 & 1.31 & 1.30 & 1.31\\
Phe & -0.11 & -0.11 & -0.11 & -0.11 & -0.11\\
Pro & -0.14 & -0.15 & -0.15 & -0.14 & -0.15\\
Ser & -0.01 & -0.02 & -0.06 & -0.01 & -0.06\\
Thr & -0.13 & -0.13 & -0.13 & -0.13 & -0.13\\
Trp & -0.03 & -0.03 & -0.03 & -0.03 & -0.03\\
Tyr & -0.09 & -0.10 & -0.11 & -0.09 & -0.10\\
Val & -0.17 & -0.17 & -0.18 & -0.17 & -0.15\\
mAb &  2e-4 & 2e-4 & 2e-4 & 2e-4 & 2e-4
	\end{tabular}
	\end{minipage}\begin{minipage}[c]{0.28\textwidth}
	\caption{Fluxes of each metabolite for medium 5, solved with varied assumed error. $w_0$, $w_{0.05}$, and $w_1$ indicate the results with $\theta$ at $0$, $5$, and $100$ \% of its given value in Table \ref{tb:alpha}, respectively.  $w_{0,inv}$ and $w_{1,inv}$ indicate the results with an interval given for \ce{CO2} with $\theta$ at $0$ and $100$ \% respectively.
		The unit of the fluxes is $\text{pmol}\cdot \text{cell}^{-1} \cdot \text{day}^{-1}$, except for Biomass for which the unit is $\text{day}^{-1}$.}
	\label{tb:M5metw}
	\end{minipage}
\end{table*}

\section{Conclusion} \label{sec:conc}
In this work we have examined the effect of errors in measurement on the solution of the EFMs-based MFA. 
The approach has been to derive a robust form of the EFMs-based MFA, a form that considers measurement errors more directly, in the sense that each value based on measurements is given an error interval and the aim is to minimize the maximum possible error on the given interval with respect to the least-squares measure. 
These types of robust optimization problems are in general not easily solved \citep{Ghaoui1997}. However, for this special case we showed that the robust form can be stated as a convex quadratic optimization problem where column generation can be applied to achieve an optimal solution. Furthermore, it was demonstrated that an external metabolite could be taken into account in the model when no measurements exists this metabolite, by considering an interval of its value.

By considering the worst-case scenario, we see how unfavourable the fit could be within the given error intervals. The case-study compares the solution of the non-robust and robust EFMs-based MFA. Those results demonstrate that the optimal solution to the robust problem is similar to the optimal solution to the non-robust formulation, \ie the EFMs-based MFA is rather robust to those assumed errors. This indicates that the errors on measurement do not induce a large change in the solution. 

Measurement errors were known to be around 20\%, therefore it was important to consider the effects of those errors, especially with respect to which EFMs are used in the optimal solution. Our second contribution, the addition of intervals, is especially relevant for metabolites for which known intervals are available in the literature but are problematic to measure in the experimental setup. In fact, the addition of intervals had a larger effect on the solution than the addition of around 20\% errors in measurements in the presented cases. 
The robust solution addressed the parameter estimation, while adding a metabolite by considering its interval, addressed the model structure. 
 We showed a way to include the knowledge of unmeasured metabolites in the column generation method. This important result allowed to achieve improvements of the model structure. This approach can be a general strategy to improve a model structure by introducing a new metabolite in a model. It can also be an approach to identify which variables should be measured when designing a new experiment.      

\renewcommand{\abstractname}{Acknowledgements}
\begin{abstract}
The work of the authors from the Department of Mathematics was supported by the Swedish Research Council. The work of the authors from the Division of Industrial Biotechnology was supported by KTH and the Swedish Governmental Agency for Innovation Systems (VINNOVA). The CHO cell line was kindly provided by Selexis (Switzerland). Culture media were kindly provided by Irvine Scientific (CA, USA). 
Finally, we thank the editor and the two anonymous referees for their valuable comments and suggestions. 
\end{abstract}
\appendix
\section{Derivation of the Robust Quadratic Program}\label{sec:DerRQP}
The robust formulation given by \eqref{eq:RobOr} can be represented as a quadratic program. To show this, consider first a reformulation of the two-norm problem,  \begin{align}
&\underset{w \geq 0}{\text{minimize}} \underset{\abs{\Delta Q_s} \leq \theta_s \abs{Q_s}}{\text{maximize}} \quad  \frac{1}{2} \norm{\mathcal{I}A_{x}E w -Q+\Delta Q} \nonumber\\
= &\underset{w \geq 0}{\text{minimize}}  \quad  \frac{1}{2}  \sum_{s=1}^n  \underset{\abs{\Delta Q_s} \leq \theta_s \abs{Q_s}}{\text{maximize}} \quad\left(\mathcal{I}A_{x}E w -Q+\Delta Q\right)_s^2,\nonumber\\
=&\underset{w \geq 0}{\text{min}}  \quad  \frac{1}{2}  \sum_{s=1}^n  \underset{\abs{\Delta Q_s} \leq \theta_s \abs{Q_s}}{\text{max}} \quad\left(\mathcal{I}A_{x}E w-Q\right)_s^2 +2\left(\mathcal{I}A_{x}E w-Q\right)_s \Delta Q_s +\Delta Q_s^2,\nonumber\\
= &\underset{w \geq 0}{\text{min}}  \quad  \frac{1}{2}  \sum_{s=1}^n \left(\mathcal{I}A_{x}E w-Q\right)_s^2 +\frac{1}{2}\underset{\abs{\Delta Q_s} \leq \theta_s \abs{Q_s}}{\text{max}}  2\left(\mathcal{I}A_{x}E w-Q\right)_s \Delta Q_s +\Delta Q_s^2. \label{eq:RBQPder1}
\end{align}
The above derivation uses that the bound on $\Delta Q$ is elementwise, thus the maximization can be moved inside the sum. 
Considering now the inner maximization, 
\begin{itemize}
\item if $(\mathcal{I}A_{x}E w -Q)_s \leq 0$ then a $\Delta Q_s^*$ that maximizes the norm is as negative as possible, \ie $\Delta Q_s^* =- \theta_s |Q_s|$,
\item if $(\mathcal{I}A_{x}E w -Q)_s \geq 0$ then a $\Delta Q_s^*$ that maximizes the norm is as positive as possible, \ie $\Delta Q_s^* = \theta_s \abs{Q_s}$.
\end{itemize}
Thus, the maximum is given by $\Delta Q_s^* = \text{sgn}((\mathcal{I}A_{x}E w -Q)_s) \theta_s |Q_s|$,
\begin{align*}
 \underset{\abs{\Delta Q_s} \leq \theta_s \abs{Q_s}}{\text{maximize}}& 2\left(\mathcal{I}A_{x}E w-Q\right)_s \Delta Q_s +\Delta Q_s^2 \\
 =&  2\left(\mathcal{I}A_{x}E w-Q\right)_s \text{sgn}((\mathcal{I}A_{x}E w -Q)_s) \theta_s \abs{Q_s} \\&+(\text{sgn}((\mathcal{I}A_{x}E w -Q)_s) \theta_s \abs{Q_s})^2,\\
 =& 2|\left(\mathcal{I}A_{x}E w-Q\right)_s  \theta_s Q_s| +(\theta_s Q_s)^2.
\end{align*}
 Then, we can set  $t_s=\abs{\left(\mathcal{I}A_{x}E w-Q\right)_s  \theta_s Q_s}$ rewriting \eqref{eq:RBQPder1} to 
 \begin{equation}\label{eq:Robwconst}
\begin{aligned}
\underset{w,t}{\text{minimize}} \quad & \frac{1}{2} \norm{\mathcal{I}A_{x}E w-Q}^2+\sum_{k=1}^n t_s+ \frac{1}{2} (\theta^TQ)^2\\
\text{subject to}\quad & t_s \geq \left(\mathcal{I}A_{x}E w-Q\right)_s  \theta_s Q_s \quad \forall s,\\
& t_s \geq -\left(\mathcal{I}A_{x}E w-Q\right)_s  \theta_s Q_s \quad \forall s,\\
&w \geq 0.
\end{aligned}
\end{equation}
The constant $(\theta^TQ)^2 $ can be disregarded from the objective function of \eqref{eq:Robwconst} and thus obtaining \eqref{eq:Robdr}.
\section{Deriving the Master- and Subproblem}\label{sc:MPSP}
In this section we will derive the master- and subproblem for the column generation of \eqref{eq:Robdr}. For the master problem, consider when the columns of $E$ are split into the known and unknown columns $E_B$ and $E_N$ respectively, the variable vector $w$ is divided in the same way, $w_B$ and $w_N$. Furthermore, $w_N$ is fixed to zero, since $E_N$ represents unknown EFMs. This gives the problem to,
\begin{equation} \label{eq:RobMPBN}
\begin{aligned}
\underset{w,t}{\text{minimize}} \quad & \frac{1}{2} \norm{\mathcal{I}A_{x}\left[E_B\; E_N\right] \begin{bmatrix}
w_B\\w_N\end{bmatrix}-Q}^2+\mathbf{1}^Tt\\
\text{subject to}\quad & t-\Theta\tilde{Q}\left(\mathcal{I}A_{x}\left[E_B\; E_N\right] \begin{bmatrix}
w_B\\w_N\end{bmatrix} -Q\right) \geq 0, \\
& t+\Theta \tilde{Q}\left(\mathcal{I}A_{x}\left[E_B\; E_N\right] \begin{bmatrix}
w_B\\w_N\end{bmatrix} -Q\right) \geq 0, \\
&w_B \geq 0,\\
&w_N=0.
\end{aligned}
\end{equation}
The formulation given by \eqref{eq:RobMPBN} is equivalent to the master problem stated in \eqref{eq:RobMP}, by inserting $w_N=0$ into the objective function and constraints in \eqref{eq:RobMPBN}, \eqref{eq:RobMP} can be obtained. 

The subproblem comes from comparing the optimality conditions of \eqref{eq:RobMPBN} and \eqref{eq:Robdr}. Because \eqref{eq:Robdr} is convex, any point that satisfies the optimality conditions of \eqref{eq:Robdr} is an optimal solution. 
Thus, the optimality conditions of \eqref{eq:Robdr} and \eqref{eq:RobMPBN} can be compared, identifying when a point that satisfies the optimality conditions of \eqref{eq:RobMPBN} also satisfies the optimality conditions for \eqref{eq:Robdr}. The optimality conditions of \eqref{eq:Robdr} are given as,
\begin{subequations}\label{eq:condrobust}
\begin{align}
 w \geq & 0, \\
 \lambda \geq &0, \label{eq:lampos} \\
 \lambda_m \geq &0, \\
  \lambda_p \geq & 0, \\
 \lambda^Tw= & 0, \\
  \lambda_{m}^T(t-\Theta \tilde{Q}\mathcal{I}A_{x}Ew +\Theta\tilde{Q}Q)=&0, \\
   \lambda_{p}^T(t+\Theta \tilde{Q}\mathcal{I}A_{x}Ew-\Theta\tilde{Q}Q)=&0, \\
\begin{bmatrix}
 E^TA_{x}^T\mathcal{I}^T(\mathcal{I}A_{x}Ew-Q)\\
 \mathbf{1}
\end{bmatrix}
- 
\begin{bmatrix}
\lambda \\
\mathbf{0}
\end{bmatrix}
-
\begin{bmatrix}
-E^TA_{x}^T\mathcal{I}^T\Theta\tilde{Q}\\
I
\end{bmatrix}
\lambda_{m}
& \nonumber \\- 
\begin{bmatrix}
E^TA_{x}^T\mathcal{I}^T\Theta\tilde{Q}\\
I
\end{bmatrix}
\lambda_{p}
=&0.
\end{align}
\end{subequations}
Note that the last condition $\lambda_m+\lambda_p=\mathbf{1}$, together with complementarity, and that the two constraints are mutually exclusive, require that exactly one constraint is active, which is consistent with what we expect the problem to do.
Consider now the optimality conditions of \eqref{eq:RobMPBN},
\begin{subequations}\label{eq:condrobustMP}
\begin{align}
 w_N =  0 \quad w_B \geq &0 , \\
\lambda_N \text{ free } \quad \lambda_B \geq &0, \\
 \lambda_m \geq &0, \\
  \lambda_p \geq & 0, \\
 \lambda_B^Tw_B= & 0, \\
  \lambda_{m}^T(t-\Theta\tilde{Q}\mathcal{I}A_{x}E_Bw_B +\Theta\tilde{Q}Q)=&0, \\
   \lambda_{p}^T(t+\Theta\tilde{Q}\mathcal{I}A_{x}E_Bw_B -\Theta\tilde{Q}Q)=&0, \\
   \mathbf{1}-\lambda_m-\lambda_p=&0,\\
\begin{bmatrix}
 E_B^TA_{x}^T\mathcal{I}^T(\mathcal{I}A_{x}E_Bw_B-Q)\\
 E_N^TA_{x}^T\mathcal{I}^T(\mathcal{I}A_{x}E_Bw_B-Q)\\
\end{bmatrix}
- 
\begin{bmatrix}
\lambda_B \\
\lambda_N
\end{bmatrix}
+
\begin{bmatrix}
E_B^TA_{x}^T\mathcal{I}^T\Theta\tilde{Q}\\
E_N^TA_{x}^T\mathcal{I}^T\Theta\tilde{Q}
\end{bmatrix}
\lambda_{m}
& \nonumber \\- 
\begin{bmatrix}
E_B^TA_{x}^T\mathcal{I}^T\Theta\tilde{Q}\\
E_N^TA_{x}^T\mathcal{I}^T\Theta\tilde{Q}\\
\end{bmatrix}
\lambda_{p}
=&0. \label{eq:consub}
\end{align}\label{eq:optcln}
\end{subequations}
The differences of these two optimality conditions come from $\lambda_N$ being free, thus, a $\lambda$ can be negative contradicting \eqref{eq:lampos}. 
Hence, the subproblem should identify if there is a $\lambda_N<0$, this can by done by an optimization problem where \eqref{eq:consub} gives an objective function. The subproblem can then be stated as to,
\begin{equation*}
 \underset{e\in E_N}{\min}e^TA_{x}^T\mathcal{I}^T(\mathcal{I}A_{x}E_Bw_B-Q +\Theta\tilde{Q}\lambda_{m}-\Theta\tilde{Q}\lambda_{p}).
\end{equation*}
The requirement that $e\in E_N$ can be fulfilled by
\begin{equation*}
\begin{aligned}
 \{e: \;A_{i}e=0, \quad
\mathbf{1}^Te \leq 1, \quad
 e_j \geq 0\; \forall j\},
\end{aligned}
\end{equation*}
 and by ensuring that the solution is an extreme point solution.
Together these constraints and objective function define an optimization problem that identifies new EFMs that can be added to the master problem \citep{Oddsdottir2014}.
\bibliographystyle{spbasic}

\end{document}